\documentclass[11pt]{article}
\usepackage{epsfig}
\usepackage{times}
\usepackage{latexsym,amssymb,amsmath,natbib,graphicx, amsthm}

\newcommand {\ME}{\mathbb{E}^{x}}

\renewcommand{\S}{\mathcal{S}}

\numberwithin{equation}{section}
\newtheorem{proposition}{Proposition}[section]
\newtheorem{corollary}{Corollary}[section]
\newtheorem{remark}{Remark}[section]
\newtheorem{lemma}{Lemma}[section]

\newcommand {\R}{\mathbb{R}}
\newcommand {\F}{\mathcal{F}}

\newcommand {\E}{\mathbb{E}}
\newcommand{\A}{\mathcal{A}}

\title{
 Optimizing Venture Capital Investments in a Jump Diffusion Model}
\author{%
Erhan Bayraktar \thanks{E. Bayraktar was supported in part by the
National Science Foundation, under grant DMS-0604491.} \and Masahiko
Egami \footnote{E. Bayraktar and M. Egami are with Department of
Mathematics, University of Michigan, Ann Arbor, MI, 48109-1043, USA.
Email:\{erhan, egami\}@umich.edu. We thank Savas Dayanik for
valuable comments and are most grateful to Anthea Au Yeung, Ikumi Koseki, T\={o}ru
Masuda, and Akihiko Yasuda for insight and
good memories through helpful and enjoyable conversations.}}
\date{}

\begin{document}

\maketitle
\begin{abstract}
\noindent We study two practical optimization problems in relation
to venture capital investments and/or Research and Development
(R\&D) investments. In the first problem, given the amount of
the initial investment and the cash flow structure at the initial
public offering (IPO), the venture capitalist wants to
maximize overall discounted cash flows after subtracting
subsequent investments, which keep the invested company solvent. We describe this problem as a
mixture of singular stochastic control and optimal stopping
problems.  The singular control corresponds to finding an optimal
subsequent investment policy so that the value of the investee
company stays solvent. The optimal stopping corresponds to finding
an optimal timing of making the company public. The second problem
is concerned with optimal dividend policy. Rather than selling the
company at an IPO, the investor may want to harvest technological
achievements in the form of dividend when it is appropriate. The
optimal control policy in this problem is a mixture of singular
and impulse controls.

\end{abstract}
\noindent \small{Key words: Venture capital investments, R\&D,
IPO, stochastic control, optimal stopping, singular control,
impulse control, jump diffusions.\\ \noindent JEL Classification:
O32, O33, C61.\\ \noindent Mathematics Subject Classification
(2000) : Primary: 49N25, Secondary: 60G40.}\\

\section{Introduction}
In accordance with the recent theoretical and practical
development in the area of real options, modeling venture capital
investment and R\&D has become increasingly an important topic,
see Davis et al. \cite{davis-2002} for a review of this literature. One of the most important issues is modeling the dynamics of
the value process of start-up companies and/or R\&D projects.
Among many approaches, one approach is to use jump models with
Poisson arrivals. For example, Willner \cite{willner} uses a
deterministic drift component and stochastic jumps whose size
follows a gamma distribution. A similar model is presented by
Pennings and Lint \cite{pennings-lint}, who also model with a
deterministic drift and a jump part whose size follows a Weibull
distribution
with scale parameter two. 
In the spirit of these papers, we will model the value of the
process with a jump diffusion. More specifically, we are assuming
that the company or the R\&D project has (unproven) innovative
technologies and hence the appreciation of the company value
occurs when there is a technological
breakthrough or discovery of innovative methods.

Let $\{\Omega,\mathcal{F},\mathbb{P}\}$ be a probability space
hosting a Poisson random measure $N(dt,dy)$ on $[0,\infty)\times
\R$ and Brownian motion $W=(W_t)_{t\ge 0}$, adapted to some
filtration $\mathbb{F}=\{\mathcal{F}_t\}_{t\ge 0}$ satisfying the
usual conditions. The mean measure of $N$ is $\nu(dt,dy)\triangleq
\lambda dt F(dy)$, where $\lambda>0$ is constant, and $F(dy)$ is
the common distribution of the jump sizes. The (uncontrolled)
value process $X^0$ of the invested company is described as follows:
\begin{equation}\label{eq:processjump}
dX_t^0=\mu(X^0_t)dt  +\sigma(X^0_t)dW_t+ \int_0^\infty y N(dt, dy), \quad t \leq \tau^{(0)},
\end{equation}
in which
\begin{equation}\label{eq:ruin-time}
\tau^{(0)}\triangleq \inf\{t \geq 0; X^{0}_t<0\}.
\end{equation}
In this set up the jumps of $X^0$ come from a compound Poisson process
whose jump size distribution is $F$. To obtain explicit results we
will take this distribution to be exponential. In
practice, investments in promising start-up companies are made
through venture capital funds (often called ``private equity fund"
as well) that raise capital from institutional investors such as
banks, insurance companies, university endowment funds, and
pension/retirement funds.  Venture capital funds screen out those
start-up companies and select several companies to invest in. Venture capitalist
allocates certain amounts of money to each promising companies to diversify risks.  As a result, for each investment project, venture capitalist has a certain initial budget.
In
many cases, the venture capital funds actively help and advise
them by taking a seat on the board of the invested companies. When
there is a technological breakthrough, this jump is materialized
in the following way:  the company and the venture capital fund
reevaluate the value of the company stock by using expected cash
flow methods provided that the company is successful in
manufacturing real products or by using comparable transactions in
the past. As a result, some new investors may become willing to
invest in the company at the re-evaluated price in a ``second
round" funding. Hence the appreciation of the stock value can be
modelled by an arrival of jumps.  The size and arrival rate of
jumps can be estimated by track records of the venture capital
funds. The final objective of these venture capital investments
is, in many cases, to make the company public through initial
public offerings (IPO's) or to sell to a third party at a premium.
However, in due course, there are times when the start-up company
faces the necessities to solicit new (additional) capital. In
turn, the venture capitalist has to make decisions on whether to
make additional investments.  We refer to this type of problem as
the ``IPO problem".

Let us mention some advantages of using a jump diffusion model
rather than a piecewise deterministic Markov model as in other works
in the literature.  As we discussed in the previous paragraph, until
going public in the IPO market, the start-up company evolves while
proving the merits and applicabilities of their technology.  At an
early stage, the company's growth mostly depends on the timing and
magnitude of jump part in (\ref{eq:processjump}).  At the time when
the company invites ``second" and ``third round" investors, it is
often the case that they have generated some cash flows from their
operation while jumps of great magnitude are not necessarily
expected. At these stages, the diffusion part of
(\ref{eq:processjump}) is becoming increasingly influential. Hence
the jump diffusion model can represent start-up companies of various
stages by appropriately modifying the parameters of the model.

To address the issue of subsequent investments in the IPO problem,
we first solve an optimal stopping problem of a \emph{reflected}
jump diffusion.  In this problem, the venture capitalist does not
allow the company's value to go below a fixed level, say $a$, with a minimal possible effort
and
attempts to find an optimal time to IPO
(Section 2.1). Next, we solve the problem in which the venture capitalist chooses the level $a$ optimally
(Section
\ref{sec:mixture}) subject to a budget constraint. In the the process of solving this problem we also solve the min-max version of it.  In mathematical terms, this problem is a mixture of local-time control (plus an impulse applied at time 0 depending on whether the start-up company's value is initially below $a$)
and optimal stopping. The local time control is
how the venture capitalist exercises controls or interventions in
terms of additional capital infusions. The optimal stopping is,
given a certain reward function at the IPO market, to find an
optimal timing of making the company public. In summary, while
making decisions with respect to additional investments, the
venture capitalist seeks to find an optimal stopping rule in order
to maximize her return, after subtracting the present value of
her intermediate investments or capital infusion.

Another problem of interest is the following: Rather than selling
outright the interest in the start-up company or R\&D investments,
the investor may want to extract values out of the company or
project in the form of dividend until the time when the value
becomes zero. This situation may be more suitable in considering
R\&D investments because one wants to harvest technological
achievements when appropriate, while one keeps the project
running. We refer to this type of problem as the ``harvesting
problem" (Section \ref{sec:harvest}) or dividend payment problem.  We prove the optimality of a threshold policy. Optimal dividend problem for L\'{e}vy processes
with negative jumps was analyzed by Avram et al. \cite{kyprianou-dividend}. Here the L\'{e}vy process we consider has positive jumps and due to this nature of the jumps one applies a mixture of impulse and singular stochastic control: When the controlled process jumps over the optimal threshold, the controller applies impulse control, when the controlled process approaches the threshold continuously, the controller reflects the controlled process, i.e., she applies singular control.

The rest of the paper is structured as follows: In section 2, we
solve the ``IPO problem", first by setting the lower threshold level
$a$ fixed and later by allowing this level vary.  In section 3, we
solve the ``harvesting problem".  Next, we construct a candidate solution  and verify the optimality of this candidate
by showing that the conditions prescribed in the
verification lemma are all satisfied.  We also provide some static
sensitivity analysis to the model parameters.  In section 4, we give our concluding remarks and compare the values of IPO and harvesting problems.

\section{The IPO Problem}\label{sec:vcsol}
\subsection{Optimal Stopping of a Reflected Diffusion}
 The
dynamics of the value of the start-up company is described as
(\ref{eq:processjump}). After making the initial investment in the
amount of $x$, the venture capitalist can make interventions in
the form of additional investments. Hence the controlled process
$X$ is written as follows:
\begin{equation}\label{eq:controlled}
    dX_t = \mu(X_t)dt + \sigma(X_t)dW_t + \int_0^\infty y N(dt, dy)
    +dZ_t, \quad X_0=x,
    \end{equation}
in which for a given $a\geq 0$,
\begin{equation}\label{eq:local-time}
Z^a_t=(a-x)1_{\{x<a\}}+L_t, \quad t\geq 0
\end{equation}
where $L_t$ is the solution of
\begin{equation}
L_t=\int_0^{t}1_{\{X_s=a\}}dL_s, \quad t \geq 0.
\end{equation}
Note that $Z=(Z_t)_{t\ge 0}$ is a continuous, non-decreasing
(except at $t=0$) $\F_t$-adapted process. In this set up, through cash infusion or additional investments, the venture capitalist aims to keep the value of the start-up above $a$ with a minimal possible effort.

The venture capitalist's purpose is to find the best
$\mathbb{F}$-stopping time to make the company public through an
IPO to maximize the present value of the discounted future
cash-flow. We will denote the set of all $\mathbb{F}$ stopping
times by $\mathcal{S}$. The discounted future cash-flow of the
venture capitalist if he applies the control $Z$, and makes an
initial public offering at time $\tau$ is
\begin{equation}\label{eq:problem}
    J^{\tau, a}(x)\triangleq\ME\left[e^{-\alpha \tau} h(X_{\tau})-
    \int_{0}^{\tau}e^{-\alpha s}dZ^a_s\right],
    \quad \tau\in\S.
\end{equation}
We use the following notation:
\begin{equation}
\int_0^{t}e^{-\alpha t}dZ_t =Z_0+\int_{(0,t)}e^{-\alpha
t}dZ_t.
\end{equation}
We will assume that $h:
\mathbb{R_+}\rightarrow \mathbb{R_+}$ has the following form:
\begin{equation}
    h(x) \triangleq
                    rx,
 \end{equation}
with $r> 1$.   The parameter $r$ is
 determined by the IPO market.

Let us discuss the rationale of our model specification.  First
of all, valuation of IPOs is itself a very challenging subject
and, to our knowledge, no complete solutions have yet been
obtained. The pricing mechanism at the IPO market is complex,
involving uncertainties with respect to the future of the newly
publicized companies. A widely observed and recommended procedure
both by academics and market practitioners is using comparable
firm ``multiples":  The subject company's operational and
financial information is compared with those of publicly-owned
comparable companies, especially with ones newly made public.  For
example, the price-earning (P/E) ratio and/or market-to-book (M/B)
ratio are multiplied by a certain number called ``multiples" (that
may vary from industry to industry) to calculate the IPO value.
These numbers of the comparable firms serve as benchmarks. Kim and
Ritter \cite{kim-ritter} found that ``P/E multiples using
forecasted earnings result in much more accurate valuations" than
using historical earnings.  In this pricing process, the role of
investment banks (they often serve underwriters as well) is
critical.  They, together with the firm, evaluate the current
operational performance, analyze the comparable firm ``multiples",
project future earnings, assess the market demand for IPO stocks
and set the timing of IPO. This procedure inherently involves
significant degree of variation on prices and introduces
``discontinuity" of the post-IPO value from the pre-IPO value,
since post-IPO value is, to a certain extent, market driven, while
pre-IPO value is mostly company specific. This justifies our
choice of reward function $h(\cdot)$. The post IPO value $h(x)=rx$ is
strictly greater than the pre IPO value $x$.

Further evidence of this discontinuity can be obtained by the
literature about initial stock returns on IPO markets. Johnston
and Madura \cite{johnston-madura} examined the initial (i.e. first
day trading) returns of IPOs (both Internet firms and non-Internet
firms) during January 1, 1996 to December 31, 2000.  They found
that the average initial returns of Internet firms IPOs was
$78.50\%$. They also reviewed the papers by Ibbotson and Jaffe
\cite{ibbotson-jaffe}, Reilly \cite{reilly} and Ritter
\cite{ritter} and tabulated those authors' findings. Ritter
\cite{ritter}, for example, found average initial returns as high
as $48.4\%$ for IPOs that occurred during 1980-1981. Thus, it is
widely observed that the IPO companies are priced at a premium and
that these premia are the most important sources of income to the
venture capitalist. This phenomenon of abnormal initial returns is
incorporated in our model with $r>1$.


The purpose of the venture
capitalist is to determine $\tau^* \in \mathcal{S}$ such that
\begin{equation}\label{eq:valuefunc}
V(x;a) \triangleq
\sup_{\tau\in\S}J^{\tau,a}(x)=J^{\tau^*,a}(x), \quad x \geq 0.
\end{equation}
if such a $\tau^*$ exists.

\subsubsection{Verification Lemma}

\begin{lemma}\label{lem:ver}
Let us assume that $\sigma(\cdot)$ is bounded. If a non-negative function $v \in
\mathcal{C}^1(\mathbb{R}_+)$  is also
twice continuously differentiable except at countably many points,
and satisfies
\begin{enumerate}
\item [(i)]  $(\A-\alpha)v(x)\leq 0$, $x \in (a,\infty)$,
\item[(ii)] $v(x) \geq h(x)$,  $x \in (a,\infty)$,
\item[(iii)] $v(x)=x-a+v(a)$, $x \in [0,a]$, and $v'(a+)=1$.
\end{enumerate}
in which the integro-differential operator $\A$ is defined by
\begin{equation}\label{eq:infinite}
\A f(x) =\mu(x)f'(x)+\frac{1}{2}\sigma^2(x)f''(x)+\lambda
\int_0^{\infty}(f(x+y)-f(y))F(dy),
\end{equation}
\begin{equation}\label{eq:verf-ge}
\text{then} \quad v(x) \geq J^{\tau,Z^a}(x), \quad \tau \in
\mathcal{S}.
\end{equation}
Moreover, if there exists a point $b(a)$ such that
\begin{enumerate}
\item [(iv)] $(\A-\alpha)v(x)=0$ and for all $x\in [a, b(a))$, $v(x)>h(x)$ for all $x  \in (a, b(a))$,
\item [(v)] $(\A-\alpha)v(x)<0$ for all $x \in (b(a),\infty)$, $v(x)=h(x)$ for all $x \in [b(a),
\infty)$,
\end{enumerate}
then $v(x)=V(x;a)$, $x \in \mathbb{R}_+$, and $\tau_{b(a)}
\triangleq \inf\{t\ge 0; X_t \ge b(a)\}$ is optimal.
\end{lemma}
\begin{proof}
Let us define $\tau(n)\triangleq\inf\{t\geq 0; X_t\geq n\}$. Let
$\tau \in \mathcal{S}$. When we apply It\^{o}'s formula to the
semimartingale $X$ (see e.g. Jacod and Shiryaev \cite{JS}), we
obtain
\begin{equation}
\begin{split}
e^{-\alpha(\tau\wedge\tau(n)\wedge\tau_0)}v(X_{\tau\wedge
\tau(n)\wedge\tau_0})&=v(x)+\int_0^{\tau\wedge\tau(n)\wedge\tau_0}
e^{-\alpha s}(\mathcal{A}-\alpha)v(X_s)ds+1_{\{x<a\}}(a-x)
\\ &+\int_0^{\tau\wedge\tau(n)\wedge\tau_0}e^{-\alpha
s}v'(X_s)dL_s+\int_0^{\tau\wedge\tau(n)\wedge\tau_0}e^{-\alpha
s}\sigma(X_s)v'(X_s)dW_s
\\&+\int_0^{\tau\wedge\tau(n)\wedge\tau_0}\int_0^{\infty}(v(X_{s-}+y)-v(X_{s-}))(N(ds,
dy)-\nu(ds,dy))
\end{split}
\end{equation}
Rearranging this equation and after taking expectations we get
\begin{equation} \label{eq:exval}
\begin{split}
v(x)&=(x-a)1_{\{x<a\}}+\ME\left[e^{-\alpha(\tau\wedge
\tau(n)\wedge \tau_0)}v(X_{\tau\wedge \tau(n)\wedge\tau_0})-\int_0^{
\tau\wedge\tau(n)\wedge\tau_0}e^{-\alpha s}dL_s\right]
\\&
+\ME\left[\int_0^{
\tau\wedge\tau(n)\wedge\tau_0}e^{-\alpha s}(1-v'(X_s))dL_s-\int_{0}^{\tau\wedge \tau(n)\wedge\tau_0}e^{-\alpha
s}(\mathcal{A}-\alpha)v(X_s)ds \right]
\\&-\ME\left[\int_0^{\tau_0\wedge \tau(n)\wedge\tau_0} e^{-\alpha
s}(v(X_{s-}+y)-v(X_{s-}))(N(ds, dy)-\nu(ds, dy))\right]
\\&-\ME \left[\int_{0}^{\tau\wedge \tau(n)\wedge\tau_0} e^{-\alpha s}
 \sigma (X_s)v'(X_s)dW_s\right].
\end{split}
\end{equation}
Since the functions $\sigma(\cdot)$, $v(\cdot)$ and $v'(\cdot)$
are bounded on the interval $[0,n]$, the expected value stochastic integral terms vanish, and since $v'(a)=1$ the expected value of the integral with respect to $L$ also vanishes. On the other hand, the expected value of the integral with respect to the Lebesgue measure is greater than zero by Assumption (i). Therefore,
\begin{equation*}
v(x)\geq(x-a)1_{\{x<a\}}+\E^{x}\left[e^{-\alpha(\tau\wedge
\tau(n)\wedge\tau_0)}v(X_{\tau\wedge
\tau(n)\wedge\tau_0})-\int_0^{
\tau\wedge\tau(n)\wedge\tau_0}e^{-\alpha s}dL_s\right].
\end{equation*}
Equation (\ref{eq:verf-ge}) follows from the bounded and monotone
convergence theorems and assumption (ii).

On the other hand, if we substitute $\tau$ for $\tau_{b(a)}$ in the
above equations and use assumptions (iv) and (v), we get that
$v(\cdot)=J^{\tau_{b(a)},a}(\cdot)$, which proves that
$v(\cdot)=V(\cdot)$ and that $\tau_{b(a)}$ is optimal, i.e.,
$J^{\tau_{b(a)},a}(\cdot) \geq J^{\tau,a}(\cdot)$ for any $\tau
\in \mathcal{S}$.
\end{proof}

\subsubsection{Construction of a Candidate Solution}\label{sec:cand1}

We will assume that the mean measure of the Poisson random measure
$N$ is given by $\nu(dt, dy)=\lambda dt \eta e^{-\eta y}dy$. In
other words, we consider the case in which the jumps come from a
compound Poisson process with exponentially distributes jump
sizes. We also assume that $\mu(x)=\mu$ where $\mu \in \R$ and
$\sigma(x)=\sigma>0$. We will also assume that $\mu+\lambda/\eta>0$. This assumption simply says that the overall trend of the company is positive which motivates the venture capitalist to keep the start-up alive.

 \noindent The action of the infinitesimal generator
of $X^0$ on a test function $f$ is given by
 \begin{equation}\label{eq:generator}
\mathcal{A}f(x)=\mu f'(x)+\frac{1}{2}\sigma^2f''(x)+
\lambda\int_0^\infty (f(x+y)-f(x))\eta e^{-\eta y}dy.
 \end{equation}
Let us define
\begin{equation}
G(\gamma) \triangleq \frac{1}{2}\sigma^2 \gamma^2+\mu \gamma+
\frac{\lambda \eta}{\eta -\gamma}-\lambda.
\end{equation}
Note that
\begin{equation}
\ME\left[e^{\gamma X_t^0}\right]=\exp\left(G(\gamma)t\right).
\end{equation}
\begin{lemma}
The equation $G(\gamma)=\alpha$ has two positive roots $\gamma_1$,
$\gamma_2$ and one negative root $-\gamma_3$ satisfying
\begin{equation}
0<\gamma_1<\eta<\gamma_2, \quad \text{and} \quad \gamma_3>0.
\end{equation}
\end{lemma}

\begin{proof}
Let us denote
\begin{equation}
A(\gamma) \triangleq \frac{1}{2} \sigma^2 \gamma^2+ \mu \gamma
-(\lambda+\alpha), \quad B(\gamma) \triangleq \frac{\lambda
\eta}{\gamma-\eta}.
\end{equation}
It follows that
\begin{equation}\label{eq:limits-B}
\lim_{\gamma \downarrow \eta}B(\gamma)=\infty, \quad \lim_{\gamma
\uparrow \eta}B(\gamma)=-\infty, \quad \lim_{\gamma \rightarrow
-\infty}B(\gamma)=0,
\end{equation}
\begin{equation}\label{eq:limits-A}
\lim_{\gamma \rightarrow -\infty}A(\gamma)=\lim_{\gamma
\rightarrow \infty}A(\gamma)=\infty,
\end{equation}
and that
\begin{equation}\label{eq:comp-at-zero}
A(0)=-(\lambda+\alpha)<B(0)=-\lambda.
\end{equation}
Moreover, $A(\cdot)$ is strictly decreasing on $(-\infty,
-\mu/\sigma^2)$ and strictly increasing on $(-\mu/\sigma^2,\infty)$;
$B(\cdot)$ is strictly decreasing both on $(-\infty,\eta)$ and on $(\eta,\infty)$ with different asymptotic behavior on different sides of $\gamma=\eta$. The claim is a direct consequence of
these observations.
\end{proof}

Let us define
\begin{equation}\label{eq:v0-defn}
 v_0(x;a) \triangleq A_1 e^{\gamma_1 x}+A_2 e^{\gamma_2 x}+ A_3 e^{-\gamma_3
x},
\end{equation}
for some $A_1, A_2, A_3 \in \mathbb{R}$ and $b>a$, which are to be
determined. We set the candidate value function as
\begin{equation}\label{eq:defn-candidate}
v(x;a) \triangleq \begin{cases} x-a+v_0(a;a), & x \in [0,a],
\\v_0(x;a),& x \in [a,b],
\\ rx, & x \in [b,\infty),
\end{cases}
\end{equation}
 Our aim is to determine
these constants so that  $v(\cdot;a)$ satisfies the conditions of
the verification lemma.

We will choose $A_1, A_2, A_3 \in \mathbb{R}$ and $b > a$ to satisfy
\begin{subequations}
\begin{eqnarray}
 &A_1e^{\gamma_1b}+A_2e^{\gamma_2b}+A_3e^{-\gamma_3b}=rb,
 \label{eq:no1}\\&
 \frac{A_1\eta}{\gamma_1-\eta}e^{\gamma_1b}+\frac{A_2\eta}{\gamma_2-\eta}e^{\gamma_2b}+\frac{A_3\eta}
{-\gamma_3-\eta}e^{-\gamma_3b} +r\left(b+\frac{1}{\eta}\right)=0,
\label{eq:no2} \\
&\gamma_1A_1e^{\gamma_1a}+\gamma_2A_2e^{\gamma_2a}-\gamma_3A_3e^{-\gamma_3a}=1,
\label{eq:no3}\\
   &\gamma_1A_1e^{\gamma_1b}+\gamma_2A_2e^{\gamma_2b}-\gamma_3A_3e^{-\gamma_3b}=r
   \label{eq:no4}.
\end{eqnarray}
\end{subequations}
For the function $v$ in (\ref{eq:defn-candidate}) to be
well-defined, we need to verify that this set of equations have a
unique solution. But before let us point how we came up with these
equations. The expressions (\ref{eq:no1}), (\ref{eq:no3}) and
(\ref{eq:no4}) come from \emph{continuous pasting} at $b$,
\emph{first-order smooth pasting} at $a$ and \emph{first order
smooth pasting at $b$}, respectively. Equation (\ref{eq:no2}) on
the other hand comes from evaluating
\begin{equation}
(\mathcal{A}-\alpha)v(x;a)=\mu
v_0'(x)+\frac{1}{2}\sigma^2v_0''(x)+\lambda\left(\int_0^{b-x}v_0(x+y)
\eta e^{-\eta y}dy+\int_{b-x}^{\infty}r\cdot(x+y) \eta e^{-\eta
y}dy \right)-(\lambda+\alpha)v_0(x)=0.
\end{equation}

\begin{lemma}\label{lem:eu-se}
For any given $a$, there is a unique $(A_1,A_2,A_3,b)\in
\mathbb{R}^{3} \times (a,\infty)$ that solves the system of
equations (\ref{eq:no1})-(\ref{eq:no4}). Moreover,
$b>\max\{a,b^*\}$, in which
\begin{equation}
b^*\triangleq
\frac{1}{\alpha}\left(\mu+\frac{\lambda}{\eta}\right)
\end{equation}
and $A_1>0$, $A_2>0$.
\end{lemma}

\begin{proof}
Using (\ref{eq:no1}), (\ref{eq:no2}) and (\ref{eq:no4}) we can
determine $A_1$, $A_2$ and $A_3$ as functions of $b$:
\begin{equation}\label{eq:exp-A123}
\begin{split}
A_1(b)&=\frac{r}{\eta^2}\frac{(\eta-\gamma_1)[\gamma_2 \gamma_3
(\eta
b+1)+\eta(\gamma_2-\gamma_3)]}{(\gamma_3+\gamma_1)(\gamma_2-\gamma_1)}\,e^{-\gamma_1
b}=:D_1(b)e^{-\gamma_1 b},
\\A_2(b)&=\frac{r}{\eta^2}\frac{(\gamma_2-\eta)[\gamma_1 \gamma_3
(\eta
b+1)+\eta(\gamma_1-\gamma_3)]}{(\gamma_3+\gamma_2)(\gamma_2-\gamma_1)}\,e^{-\gamma_2
b}=:D_2(b)e^{-\gamma_2 b},
\\A_3(b)&=\frac{r}{\eta^2}\frac{(\eta+\gamma_3)[\gamma_2 \gamma_1
(\eta
b+1)-\eta(\gamma_2+\gamma_1)]}{(\gamma_3+\gamma_1)(\gamma_2+\gamma_3)}\,e^{\gamma_3
b}=:D_3(b)e^{\gamma_3 b}.
\end{split}
\end{equation}
Let us define
\begin{equation}
R(b) \triangleq \gamma_1 A_1(b)e^{\gamma_1 a}+\gamma_2
A_2(b)e^{\gamma_2 a}-\gamma_3 A_3(b)e^{-\gamma_3 a}.
\end{equation}
To verify our claim we only need to show that there is one and only
one root of the equation $R(b)=1$. 
Observe that
\begin{equation}
\begin{split}
R(a)&=\gamma_1 A_1(a)e^{\gamma_1 a}+\gamma_2 A_2(a)e^{\gamma_2
a}-\gamma_3 A_3(a)e^{-\gamma_3 a}=\gamma_1 D_1(a)+\gamma_2
D_2(a)-\gamma_3 D_3(a)=r>1,
\end{split}
\end{equation}
and that
\begin{equation}
\lim_{b \rightarrow \infty}R(b)=-\infty.
\end{equation}

The derivative of $b \rightarrow R(b)$ is
\begin{equation}\label{eq:derivative-R}
\begin{split}
R'(b)&=\gamma_1 A'_1(b)e^{\gamma_1 a}+\gamma_2 A'_2(b)e^{\gamma_2
a}-\gamma_3 A'_3(b)e^{-\gamma_3 a}
\\&=\left[\gamma_1 C_1 e^{-\gamma_1(b-a)}+\gamma_2 C_2
e^{-\gamma_2(b-a)}+\gamma_3 C_3 e^{\gamma_3(b-a)}\right](-\eta
\gamma_1 \gamma_2 \gamma_3 b+Y),
\end{split}
\end{equation}
in which
\begin{equation}
C_1 \triangleq
\frac{r}{\eta^2}\frac{\eta-\gamma_1}{(\gamma_3+\gamma_1)(\gamma_2-\gamma_1)}>0,
\,\, C_2 \triangleq
\frac{r}{\eta^2}\frac{\gamma_2-\eta}{(\gamma_3+\gamma_2)(\gamma_2-\gamma_1)}>0,
\,\, \text{and}\,\, C_3 \triangleq
\frac{r}{\eta^2}\frac{\eta+\gamma_3}{(\gamma_3+\gamma_1)(\gamma_3+\gamma_2)}>0,
\end{equation}
\begin{equation}
\text{and} \quad Y \triangleq -\gamma_1 \gamma_2 \gamma_3+\eta
(-\gamma_1 \gamma_2+\gamma_2 \gamma_3+ \gamma_1 \gamma_3).
\end{equation}
Observe that
\begin{equation}
 \frac{Y}{\eta \gamma_1 \gamma_2
\gamma_3}=b^*.
\end{equation} From (\ref{eq:derivative-R}) it
follows that on $(-\infty,b^*]$ the function $b \rightarrow R(b)$
is increasing, and on $[b^*,\infty)$ it is decreasing. If $b^*
\leq a$, then it follows directly from  $R(a)=r>1$ and $\lim_{b
\rightarrow \infty}R(b)=-\infty$ that there exists a unique $b>a$
such that $R(b)=1$. On the other hand, if $b^*>a$, then $R(x)>1$
on $x \in [a,b^*]$. Again, since $\lim_{b \rightarrow
\infty}R(b)=-\infty$, there exists a unique $b>b^*$ such that
$R(b)=1$.

Let us show that $A_1(b)>0$ for the unique root of $R(b)=1$.
Observe that $A_1'(b^*)=0$ and $A_1(b^*)>0$. Moreover, $b^*$ is
the only local extremum of the function $b \rightarrow A_1(b)$,
and $\lim_{b \rightarrow \infty}A_1(b)=0$. Since this function is
decreasing on $[b^*,\infty)$, $A_1(b)>0$. Similarly, $A_2(b)>0$.
\end{proof}

\begin{remark}\label{rem:smft}
 It
follows from (\ref{eq:no1}),(\ref{eq:no3}) and (\ref{eq:no4}) that
\begin{equation}\label{eq:inc-a-b}
v(b;a)=rb, \quad v'(a;a)=1<v'(b(a);a)=r.
\end{equation}
\end{remark}

\begin{lemma}\label{lem:convexi}
Let $A_1$, $A_2$, $A_3$, and $b$ be as in Lemma~\ref{lem:eu-se}
and $v_0(\cdot;a)$ be as in (\ref{eq:v0-defn}). Then if $A_3 \geq
0$, then $v_0(\cdot;a)$ is convex for all $a \geq 0$. Otherwise,
there exists a unique point $\tilde{x}<b$ such that,
$v_0(\cdot;a)$ is concave on $[0,\tilde{x})$ and convex on
$(\tilde{x},\infty)$.
\end{lemma}

\begin{proof}
The first and the second derivative of $v_0(\cdot;a)$ (defined in
(\ref{eq:defn-candidate})) are
\begin{equation}\label{eq:derivatives}
\begin{split}
v_0'(x;a)&=A_1 \gamma_1 e^{\gamma_1 x}+A_2 \gamma_2 e^{\gamma_2
x}- \gamma_3 A_3 e^{-\gamma_3 x}, \quad  v_0''(x;a)=A_1 \gamma_1^2
e^{\gamma_1 x}+A_2 \gamma_2^2 e^{\gamma_2 x}+ \gamma_3^2 A_3
e^{-\gamma_3 x}.
\end{split}
\end{equation}

>From Lemma~\ref{lem:eu-se} we have that
\begin{equation}
A_1>0 \quad \text{and} \quad A_2>0.
\end{equation}
If $A_3 \geq 0$, then (\ref{eq:derivatives}) and
(\ref{eq:inc-a-b}) imply that $v''(x;a)>0$, $x \in [a,b(a)]$,
i.e., $v(\cdot;a)$ is convex on $[a,b(a)]$.

Let us analyze the case when $A_3<0$. In this case the functions
$x \rightarrow A_1 \gamma_1^2 e^{\gamma_1 x}+A_2 \gamma_2^2
e^{\gamma_2 x}$ and $x \rightarrow -\gamma_3^2 A_3 e^{-\gamma_3
x}$ intersect at a unique point $\tilde{x}>0$. The function
$v_0'(\cdot;a)$ (defined in (\ref{eq:defn-candidate})) decreases
on $[0,\tilde{x})$ and increases on $[\tilde{x},\infty)$. Now from
(\ref{eq:inc-a-b}) it follows that $\tilde{x}<b(a)$.

\end{proof}

\subsubsection{Verification of Optimality}

\begin{proposition}\label{prop:main1}
Let us denote the unique $b$ in Lemma~\ref{lem:eu-se} by $b(a)$ to
emphasize its dependence on $a$. Then
$v(\cdot;a)$ defined in (\ref{eq:defn-candidate}) is equal to
$V(\cdot;a)$ of (\ref{eq:valuefunc}).
\end{proposition}

\begin{proof}
The function $v$ in (\ref{eq:defn-candidate}) already satisfies
\begin{equation}
(\mathcal{A}-\alpha)v(x;a)=0, \quad x \in (a,b(a)), \quad
v(x;a)=rx, \quad x \in  [b(a),\infty), \quad v'(a;a)=1.
\end{equation}
Therefore, we only need to show that
\begin{equation}\label{eq:in-eqs}
(\mathcal{A}-\alpha)v(x;a)<0, \quad x \in (b(a),\infty), \quad
\text{and that} \quad  v(x;a)>rx, \quad x \in (a,b(a)),
\end{equation}
Let us prove the first inequality.
\begin{equation}
(\mathcal{A}-\alpha)v(x;a)=\mu r+ \frac{\lambda r}{\eta}-\alpha r
x, \quad x>b(a).
\end{equation}
So, $(\mathcal{A}-\alpha)v(x;a)<0$, for $x>b(a)$ if and only if
\begin{equation}\label{eq:bgbs}
b(a)>\frac{1}{\alpha}\left(\mu+\frac{\lambda}{\eta}\right)=b^*.
\end{equation}
However, we already know from Lemma~\ref{lem:eu-se} that
(\ref{eq:bgbs}) holds.

Let us prove the second inequality in (\ref{eq:in-eqs}). If $A_3 \geq 0$, then Lemma~\ref{lem:convexi} imply that
$v''(x;a)>0$, $x \in [a,b(a)]$, i.e., $v(\cdot;a)$ is convex on
$[a,b(a)]$. Therefore $v'(\cdot;a)$ is increasing on $[a,b(a)]$
and $v'(x;a) \in [1,r)$ on $[a,b)$. Since $x \rightarrow v(x;a)$
intersects the function $x \rightarrow rx$ at $b(a)$, $v(x;a)>rx$,
$x \in [a,b(a))$. Otherwise there would exist a point $x^* \in
[a,b)$ such that $v'(x^*;a)>r$.

If $A_3<0$, then the function $v_0'(\cdot;a)$ (defined in
(\ref{eq:defn-candidate})) decreases on $[0,\tilde{x})$ and
increases on $[\tilde{x},\infty)$, in which $\tilde{x}<b(a)$, by
Lemma~\ref{lem:convexi}. If $\tilde{x} \geq a$, then $v'(x;a)<r $
for $x \in [a,b(a))$ since $v'(a;a)=1$, $v'(x;a)<1$ for $x \in
(a,\tilde{x}]$ and $v'(x;a)<r$ for $x \in (\tilde{x},b(a))$. On
the other hand if $\tilde{x}<a$, then $v'(x;a) \in [1,r)$ for $x
\in [a,b(a))$ since $v'(\cdot;a)$ is increasing on this interval
and $v'(b;a)=r$. So in any case $v'(x;a)<r$ on $[a,b(a))$. Since
$v(b;a)=rb$, then $v(x)>rx$, $x \in [a,b)$. Otherwise  there would
exist a point $x^* \in [a,b)$ such that $v'(x^*;a)>r$.
\end{proof}

\begin{figure}[h]
\begin{center}
\begin{minipage}{0.45\textwidth}
\centering \includegraphics[scale=0.6]{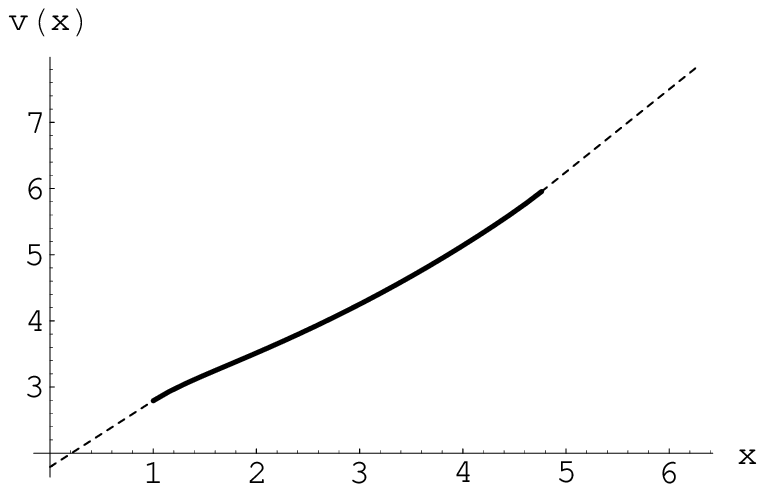} \\ (a)
\end{minipage}
\begin{minipage}{0.45\textwidth}
\centering \includegraphics[scale=0.6]{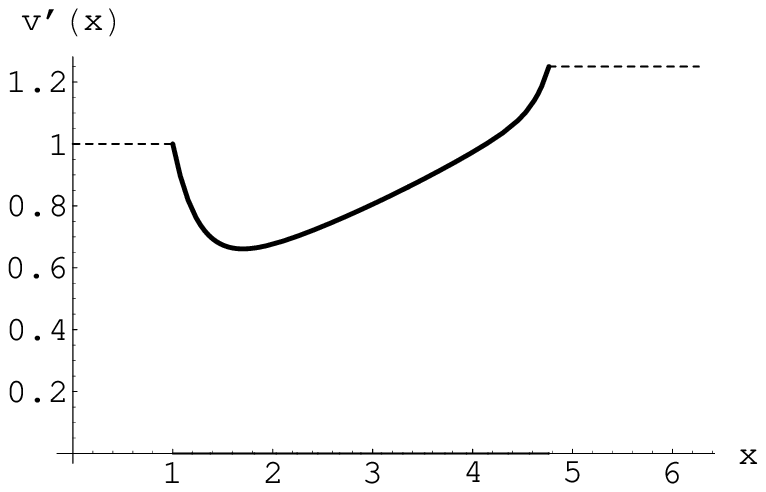} \\ (b)
\end{minipage}

\caption{\small The IPO problem with parameters
 $(\mu, \lambda, \eta, \sigma, \alpha, r)=(-0.05, 0.75, 1.5, 0.25, 0.1, 1.25)$ and $a=1$:  (a)
 The value function $v(x)$ with  $b(a)=4.7641$.
 (b) $v'(x)$ is also continuous in $x\in \R_+$.}
\label{fig:ipo}
\end{center}
\end{figure}

Figure \ref{fig:ipo} shows the value function and its derivatives. As expected the value function $v(\cdot;1)$ is concave at first and becomes convex before it coincides with the line $h(\cdot)$. It can be seen that $v(\cdot;a)$ satisfies the conditions of the verification lemma.

\subsection{Maximizing Over the Cash-Infusion Level $a$}\label{sec:mixture}

In this section, the goal of the venture
 capitalist is to find an $a^* \in [0,B]$ and $\tau^* \in \mathcal{S}$
 such that
\begin{equation}\label{eq:valuefunc-m}
U(x) \triangleq \sup_{a \in [0,B]
}\sup_{\tau\in\S}J^{\tau,a}(x)=\sup_{a \in [0, B]}V(x;a)=
J^{\tau^*,a^*}(x), \quad x \geq 0,
\end{equation}
if $(a^*,\tau^*) \in [0,B] \times \mathcal{S}$ exists.
In this optimization problem, the constraint $a
\leq B$, reflects the fact that the venture capitalist has a finite initial budget to pump-up the value
of the start-up company: the first term in (\ref{eq:local-time}) can not be greater than $B$.
The main result of this section is
Proposition~\ref{prop:mossc}. We will show that $V(x;a)$, for all
$x \geq 0$, is maximized at either $a=0$ or $a=B$. In the mean time we will also find a solution to the min-max problem
\begin{equation}\label{eq:tilde-U}
\tilde{U}(x) \triangleq \inf_{a \in [0,B]
}\sup_{\tau\in\S}J^{\tau,a}(x).
\end{equation}

We will start with analyzing the local extremums of the function
$a \rightarrow v_0(x;a)$, $x \geq 0$. We will derive the second
order smooth fit condition at $a$, from a first order derivative
condition.

\begin{lemma}\label{lem:nec-cond}
Recall the definition of the function $v_0(\cdot;a)$ from
(\ref{eq:defn-candidate}). If $\tilde{a} \geq 0$ is a local
extremum of the function $a \rightarrow v_0(x;a)$, for any $x \geq
0$, then $v_0''(\tilde{a};\tilde{a})=0$.
\end{lemma}

\begin{proof}
Let us denote
\begin{equation}\label{eq:defn-tildeA-s}
\tilde{A}_1(a) \triangleq A_1(b(a)), \quad \tilde{A}_2(a)
\triangleq A_2(b(a)), \quad \text{and} \quad \tilde{A}_3(a)
\triangleq A_3(b(a)),
\end{equation}
in which the functions $A_1(\cdot)$, $A_2(\cdot)$ and $A_3(\cdot)$
are given by (\ref{eq:exp-A123}). The derivative
\begin{equation}\label{eq:der-v0}
\frac{d v_0}{da}(x;\tilde{a})=\tilde{A}'_1(\tilde{a}) e^{\gamma_1
x}+\tilde{A}'_2(\tilde{a}) e^{\gamma_2 x}+\tilde{A}'_3(\tilde{a})
e^{-\gamma_3 x}=0
\end{equation}
for all $x \geq 0$ if and only if
\begin{equation}\label{eq:der-A-til}
\tilde{A}_1'(\tilde{a})=\tilde{A}_2'(\tilde{a})=\tilde{A}_3'(\tilde{a})=0,
\end{equation}
since the functions $x \rightarrow e^{\gamma_1 x}$, $x \rightarrow
e^{\gamma_2 x}$ and $x \rightarrow e^{-\gamma_2 x}$, $x \geq 0$,
are linearly independent.

It follows from (\ref{eq:no3}) that for any $a \geq 0$
\begin{equation}
\gamma_1 \tilde{A}_1(a)e^{\gamma_1a}+\gamma_2
\tilde{A}_2(a)e^{\gamma_2a}-\gamma_3\tilde{A}_3(a)e^{-\gamma_3a}=1.
\end{equation}
Taking the derivative with respect to $a$ we get
\begin{equation}
(\gamma_1^2 \tilde{A}_1(a)+\gamma_1 \tilde{A}'_1(a))
e^{\gamma_1a}+(\gamma_2^2
\tilde{A}_2(a)+\tilde{A}'_2(a))e^{\gamma_2a}+(\gamma_3^2\tilde{A}_3(a)-\gamma_3
\tilde{A}_3'(a)) e^{-\gamma_3a}=0.
\end{equation}
Evaluating this last expression at $a=\tilde{a}$ we obtain
\begin{equation}
\gamma_1^2 \tilde{A}_1(\tilde{a}) e^{\gamma_1
\tilde{a}}+\gamma_2^2 \tilde{A}_2(\tilde{a})e^{\gamma_2 \tilde{
a}}+\gamma_3^2\tilde{A}_3(\tilde{a}) e^{-\gamma_3
\tilde{a}}=v_0''(\tilde{a};\tilde{a})=0,
\end{equation}
where we used (\ref{eq:der-A-til}).
\end{proof}

\begin{lemma}\label{lem:prop-b}
Let $\tilde{a}$ be as in Lemma~\ref{lem:nec-cond} and $a
\rightarrow b(a)$, $a \geq 0$, be as in
Proposition~\ref{prop:main1}. Then $b'(\tilde{a})=0$. The point
$\tilde{a}$ is a unique local extremum of $a \rightarrow
v_0(x;a)$, for all $x \geq 0$, if and only if $\tilde{a}$ is the
unique local extremum of $a \rightarrow b(a)$. If $\tilde{a}$ is
the unique local extremum of $b(\cdot)$, then $b''(\tilde{a})>0$.
Moreover, $\tilde{a}=\text{argmin}_{a \geq 0 }(b(a))$.
\end{lemma}

\begin{proof}
Let $\tilde{A}_1(\cdot)$ be as in (\ref{eq:defn-tildeA-s}). Since
\begin{equation}
\tilde{A}_1'(\tilde{a})=\frac{d}{d
b}A_1(b(\tilde{a}))b'(\tilde{a})=0,
\end{equation}
and for any $a$, $b(a)>b^*$, in which $b^*$ is the unique local
extremum of the function $b \rightarrow A_1(b)$, it follows that
$b'(\tilde{a})=0$.

Assume that $\tilde{a}$ is the unique local extremum of
$b(\cdot)$. Then
\begin{equation}
\tilde{A}'_1(a)=\tilde{A}'_2(a)=\tilde{A}'_3(a)=0,
\end{equation}
if and only if $a=\tilde{a}$. Using (\ref{eq:der-v0}), it is
readily seen that $a \rightarrow v_0(x;a)$ has a unique local
extremum and that this local extremum is equal to $\tilde{a}$.

On the other hand we know from Lemma~\ref{lem:eu-se} that $b(a)>a$
for all $a$. Therefore, if $b(\cdot)$ has a unique local extremum
at $\tilde{a}$, it can not be a local maximum. On the other hand,
if there were an $a \neq \tilde{a}$ such that $b(a) \leq
b(\tilde{a})$, then there would be a local maximum in
$(\min\{a,\tilde{a}\},\max\{a,\tilde{a}\})$, which yields a
contradiction.
\end{proof}

\begin{lemma}\label{lem:pm}
 Recall the definition of $v(\cdot;a)$, $a \geq 0$, from
(\ref{eq:defn-candidate}).  For any $a_1, a_2 \geq 0$, if
$b(a_1)> b(a_2)$, then $v(x;a_1) > v(x;a_2)$, $x \geq 0$.
\end{lemma}
\begin{proof}
We will first show that
\begin{equation}\label{eq:v0-tv0}
v_0(x;a_1)=\tilde{A}_1(a_1)e^{\gamma_1
x}+\tilde{A}_2(a_1)e^{\gamma_2 x}+\tilde{A}_3(a_1)e^{-\gamma_3 x}
\geq v_0(x;a_2)=\tilde{A}_1(a_2)e^{\gamma_1
x}+\tilde{A}_2(a_2)e^{\gamma_2 x}+\tilde{A}_3(a_2)e^{-\gamma_3 x},
\end{equation}
for $x \in [0,b(a_2)]$, in which $\tilde{A}_1$, $\tilde{A}_2$,
$\tilde{A}_3$ are defined in (\ref{eq:defn-tildeA-s}). Since
$b(a_2)\leq b(a_1)$,
\begin{equation}
A_1(b(a_1))<A_1(b(a_2)), \quad A_2(b(a_1))<A_2(b(a_2)), \quad
A_3(b(a_1))>A_3(b(a_2)).
\end{equation}
This follows from the fact that the functions $A_1(\cdot)$,
$A_2(\cdot)$ are increasing and $A_3(\cdot)$ is decreasing on
$[b^*,\infty)$ and that $b(a)>b^*$, for any $a \geq 0$. See
(\ref{eq:exp-A123}) and Lemma~\ref{lem:eu-se}.

Let us define
\begin{equation}
W(x) \triangleq v_0(x;a_1)-v_0(x;a_2), \quad x \in \mathbb{R}.
\end{equation}
The derivative
\begin{equation}\label{eq:der-ltz}
W'(x)<0, \quad x \in \mathbb{R}, \quad \text{and}
\end{equation}
\begin{equation}
\lim_{x\rightarrow -\infty}W(x)=\infty, \quad \lim_{x\rightarrow
\infty}W(x)=-\infty.
\end{equation}
Therefore, $W(\cdot)$ has a unique root. We will show that this
root, which we will denote by $k$, satisfies $b(a_2)<k<b(a_1)$.

It follows from Lemma~\ref{lem:convexi} that $v_0(\cdot;a)$ is
convex on $[b(a),\infty)$, for any $a$. Moreover, for any $a \geq
0$, $v_0(\cdot;a)$ smoothly touches the function $h(\cdot)$ (see
(\ref{eq:no1}) and (\ref{eq:no4})), and stays above $h(\cdot)$ since
$v_0(\cdot;a)$ is convex. Now, since $b(a_2)<b(a_1)$, for the
function $W(\cdot)$ to have a unique root, that unique root has to
satisfy $b(a_2)<k<b(a_1)$. This proves (\ref{eq:v0-tv0}).

>From (\ref{eq:v0-tv0}) it follows that $v(x;a_1) \geq v(x;a_2)$, for
any $x \in [\min\{a_1,a_2\},b(a_2)]$. But $v(x;a_2)=rx$ for $x \geq
b(a_2)$ and $v(x;a)=v_0(x;a_1)>rx$, $x \in [b(a_2),b(a_1)]$ and
$v(x;a)=rx$, $x \geq b(a_1)$. Therefore, we have
\begin{equation}\label{eq:ineq-v}
v(x;a_1) \geq v(x;a_2), \quad x \geq \min\{a_1,a_2\}.
\end{equation}
In what follows we will show that the inequality in
(\ref{eq:ineq-v}) also holds on $x \leq \min\{a_1, a_2\}$.

Let us assume that $a_1<a_2$. It follows from (\ref{eq:der-ltz})
and $v_0'(a_2;a_2)=1$ (see (\ref{eq:no3})) that
\begin{equation}\label{eq:v0a2a1l1}
v'_0(a_2;a_1)<1.
\end{equation}
Therefore, $v_0(\cdot;a_1)$ does not intersect $x \rightarrow
x-a_2+v_0(a_2,a_2)$ $x \in [a_1,a_2]$. Otherwise, at the point of
intersection, say $\hat{x}$, $v_0'(\hat{x},a_1)>1$, which together
with (\ref{eq:v0a2a1l1}) contradicts Lemma~\ref{lem:convexi}. This
implies that $v(x;a_1)> v(x;a_2)$, $x \in [a_1,a_2]$.

Let us assume that $a_1>a_2$ and that $b(a_2)\geq a_1$. Then,
$v_0(x;a_2)<x-a_1+v_0(a_1;a_1)$. Otherwise, $v_0(\cdot;a_2)$
intersects $x \rightarrow x-a_1+v_0(a_1;a_1)$ at $x_0 \in
(a_2,a_1)$. Then $v_0'(x_0;a_2)>1$. Since $v_0'(a;a)=1$ for any
$a\ge 0$, using Lemma~\ref{lem:convexi}, it follows that
$v_0(a_1;a_2)>v_0(a_1;a_1)$. This yields a contradiction since
$v_0(\cdot;a_2)$ and $v_0(\cdot;a_1)$ do not intersect for any $x <
k$, in which $k>b(a_2) \geq a_1$. Therefore, $v(x;a_1)> v(x;a_2)$,
$x \geq [a_2,a_1]$.

Finally, let us assume that $a_1>a_2$ and that $b(a_2)<a_1$. Since
$v(x;a_2)=rx$ and $v_0(x;a_1)>rx$ for $x \geq [b(a_2),a_1]$ it
follows that $v(x;a_1)> v(x;a_2)$, $x \geq [a_2,a_1]$. Now, the
proof is complete.
\end{proof}

\begin{corollary}\label{lem:vtaleva}
 Recall the definition of $v(\cdot;a)$ from
(\ref{eq:defn-candidate}). Let $\tilde{a}$ be the unique local
extremum of $a \rightarrow b(a)$, $a \geq 0$. Then $v(x;\tilde{a})
\leq v(x;a)$, $x \geq 0$, for all $a \geq 0$.
\end{corollary}

\begin{proof}
The proof follows from Lemmas~\ref{lem:prop-b} and \ref{lem:pm}.
\end{proof}

\begin{corollary}
 Let $\tilde{a}$ be the unique local
extremum of $a \rightarrow b(a)$, $a \geq 0$. Then function
$v(\cdot;\tilde{a})$ is convex. Moreover, $v'(x;\tilde{a})>1$,
$x>\tilde{a}$.
\end{corollary}
\begin{proof}
Let $\tilde{A}_3(\tilde{a})$ be as in (\ref{eq:defn-tildeA-s}). If
$\tilde{A}_3(\tilde{a}) \geq 0$ then $v_0(\cdot;\tilde{a})$ is
convex by Lemma~\ref{lem:convexi}.

If $\tilde{A_3}(\tilde{a})<0$, then the function
\begin{equation}\label{eq:td}
v_0'''(x;\tilde{a})=\tilde{A}_1(\tilde{a}) \gamma_1^3 e^{\gamma_1
x}+\tilde{A}_2(\tilde{a}) \gamma_2^3 e^{\gamma_2 x}- \gamma_3^3
\tilde{A}_3(\tilde{a}) e^{-\gamma_3 x}>0.
\end{equation}
Since $v_0''(\tilde{a},\tilde{a})=0$, then (\ref{eq:td}) implies
that $v_0''(x;\tilde{a})>0$ for $x>\tilde{a}$. The convexity of
$v(\cdot;\tilde{a})$ follows, since it is equal to
$v_0(\cdot;\tilde{a})$ on $[\tilde{a},b(\tilde{a})]$ and is linear
everywhere else.

Since $v'(\tilde{a};\tilde{a})=1$ (see Remark~\ref{rem:smft}), it
follows from the convexity of $v(\cdot;\tilde{a})$ that
$v'(x;\tilde{a})>1$ for $x>\tilde{a}$.

\end{proof}

Note that the second order smooth fit condition $v''(\tilde{a};\tilde{a})=0$ yields a solution that minimizes $V(x;a)$, $x \geq 0$, $a \geq 0$, as a result of Corollary~\ref{lem:vtaleva}. In the next proposition we find the maximizer.

\begin{proposition}\label{prop:mossc}
Assume that  $a \rightarrow b(a)$, $a \geq 0$ has a unique local
extremum at $\tilde{a}$. Then
\begin{equation}
U(x)=\max_{a \in \{0,B\}}v(x;a), \quad \text{and} \quad \tilde{U}(x)=v(x;\tilde{a}),
\end{equation}
in which $U$ and $\tilde{U}$ are given by (\ref{eq:valuefunc-m}) and (\ref{eq:tilde-U}), respectively.
\end{proposition}
\begin{proof}
It follows from Lemma~\ref{lem:prop-b} that $a \rightarrow b(a)$,
 has a unique local extremum, and in fact this local
extremum is a minimum. Therefore, $a \rightarrow b(a)$, $a \in
[0,B]$ is maximized at either of the boundaries. The result
follows from Lemma~\ref{lem:pm}.
\end{proof}

\begin{figure}[h]
\begin{center}
\begin{minipage}{0.3\textwidth}
\centering \includegraphics[scale=0.6]{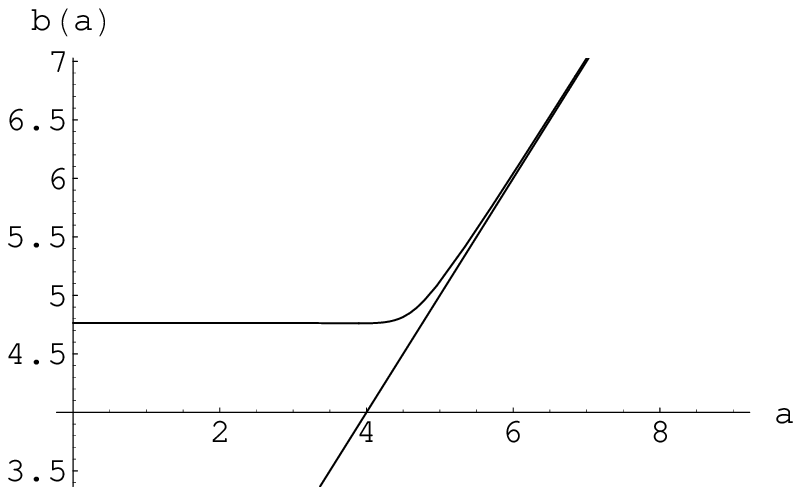} \\ (a)
\end{minipage}
\begin{minipage}{0.3\textwidth}
\centering \includegraphics[scale=0.6]{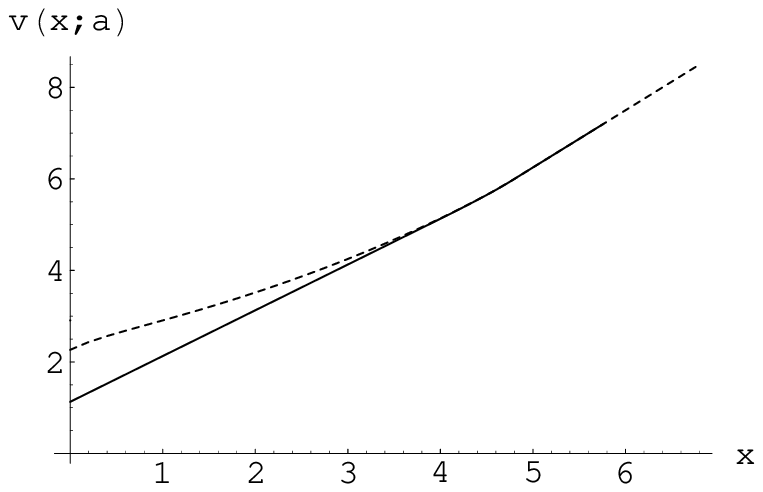} \\ (b)
\end{minipage}
\begin{minipage}{0.3\textwidth}
\centering \includegraphics[scale=0.65]{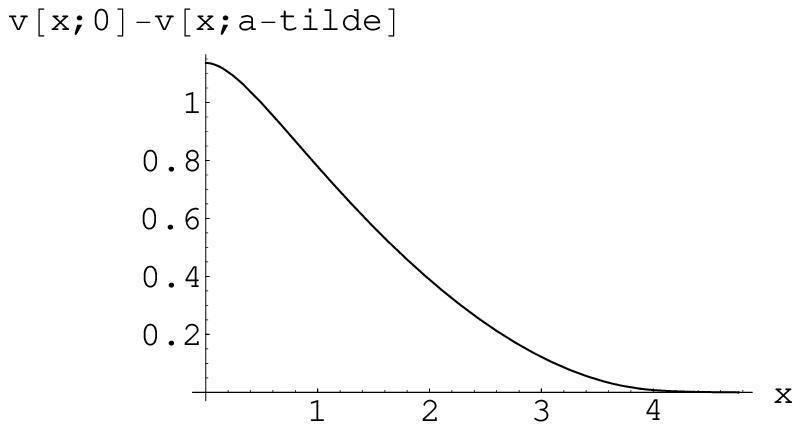} \\ (c)
\end{minipage}
\caption{\small Using the parameters $(\mu, \lambda, \eta, \sigma,
\alpha)=(-0.05, 0.75, 1.5, 0.25, 0.1)$:  (a) $\tilde{a}=3.884$
minimizes the function $b(a)$ with $b(\tilde{a})=4.741$.  (b) The
corresponding value function $v(x; \tilde{a})$ (solid line) is below
$v(x; 0)$ (dashed line). (c) $v(x;0)-v(x;\tilde{a})$.}
 \label{fig:minimum-a-tilde}
\end{center}
\end{figure}

Using the same parameters as in Figure \ref{fig:ipo}, we solve
(\ref{eq:no1})-(\ref{eq:no4}) and
\begin{equation}\label{eq:conjecture}
v_0''(\tilde{a};\tilde{a})=\gamma_1^2 \tilde{A}_1(\tilde{a})
e^{\gamma_1 \tilde{a}}+\gamma_2^2
\tilde{A}_2(\tilde{a})e^{\gamma_2 \tilde{
a}}+\gamma_3^2\tilde{A}_3(\tilde{a}) e^{-\gamma_3 \tilde{a}}=0.
\end{equation}
numerically and find $\tilde{a}$ and confirm its uniqueness. We observe in Figure
\ref{fig:minimum-a-tilde} that (a) $\tilde{a}$ is the minimizer of
$b(a)$, and (b) $v(x; 0)\ge v(x; \tilde{a})$ for $x\in \R_+$.

\begin{figure}[h]
\begin{center}
\begin{minipage}{0.40\textwidth}
\centering \includegraphics[scale=0.65]{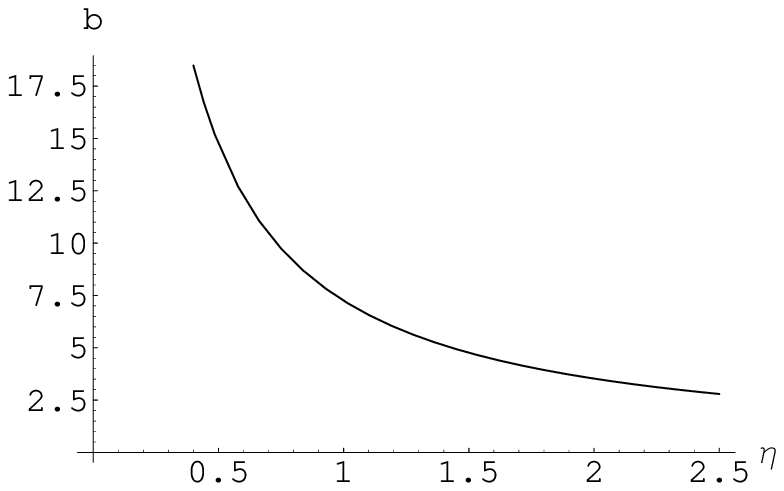} \\ (a)
\end{minipage}
\begin{minipage}{0.40\textwidth}
\centering \includegraphics[scale=0.65]{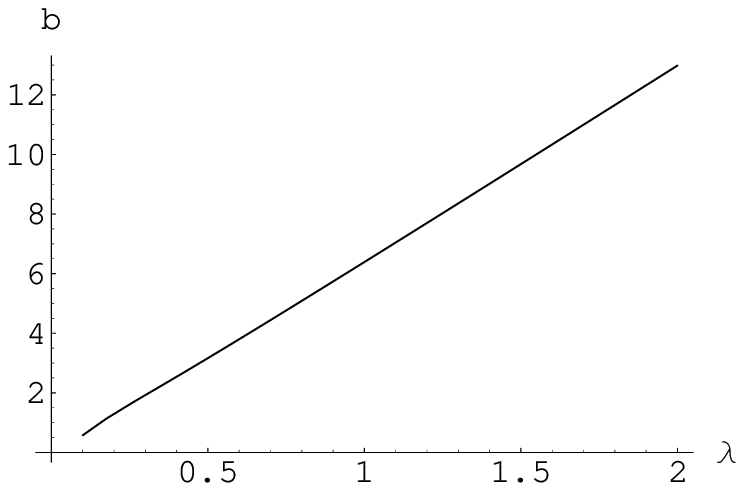} \\ (b)
\end{minipage}
\begin{minipage}{0.40\textwidth}
\centering \includegraphics[scale=0.65]{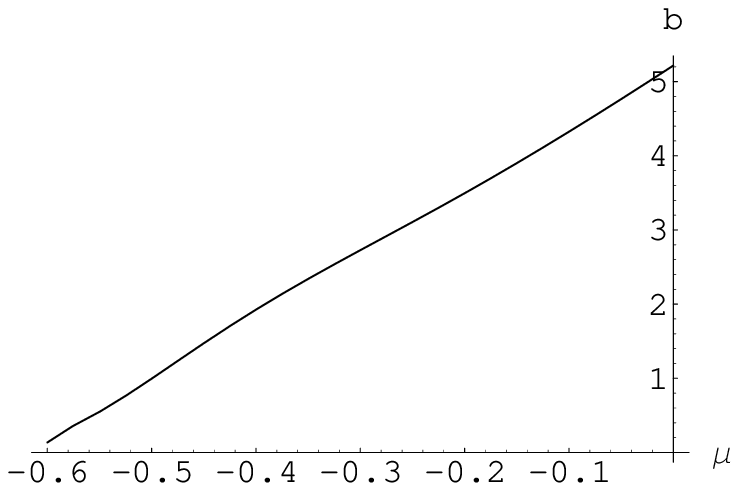} \\ (c)
\end{minipage}
\begin{minipage}{0.40\textwidth}
\centering \includegraphics[scale=0.65]{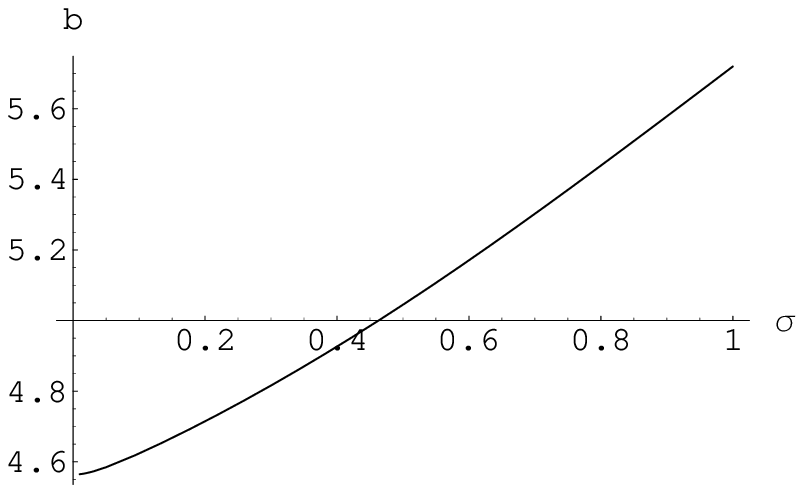} \\ (d)
\end{minipage}
\caption{\small Sensitivity analysis of the harvesting (dividend payout) problem to the parameters.  The basis parameters are $(\mu, \lambda, \eta, \sigma, \alpha)=(-0.05, 0.75, 1.5, 0.25, 0.1)$:  (a) jump  size parameter $\eta$,
(b) arrival rate $\lambda$, (c) drift rate $\mu(x)=-\mu$ and (d) volatility $\sigma$.}
 \label{fig:ipo-sensitivity}
\end{center}
\end{figure}

Before ending this section, we provide sensitivity analysis of the
optimal stopping barrier to the parameters of the problem.  We use
the parameter sets $(\mu, \lambda, \eta, \sigma, \alpha)=(-0.05,
0.75, 1.5, 0.25, 0.1)$ with $r=1.25$ and $a=0$ and vary one
parameter with the others fixed at the base case.  Figure
\ref{fig:ipo-sensitivity} shows the results. In fact, all the graphs
show monotone relationship between $b(a)$ and the parameters, which
is intuitive.  Larger $\eta$ (that means smaller $1/\eta$) leads to
a smaller threshold value since the mean jump size is small (Graph
(a)).  Similarly, larger $\lambda$ leads to a larger threshold value
since the frequency of jumps is greater and the investor can expects
higher revenue.  (Graph (b)).  In the same token, if the absolute
value of $\mu$ is greater (when the drift is negative), the process inclines to return to zero
more frequently.  Hence the investor cannot expect high revenue due
to the time value of money.  (Graph (c)).  A larger volatility expands the continuation region since the  process $X$ has a greater probability to reach further out within a fixed amount of time. Hence the investor can expect the process to reach a higher return level (Graph (d)).

\section{The Harvesting Problem} \label{sec:harvest}
\subsection{Problem Description}
In this section, the investor wants to extract the value out of
the company intermittently (i.e., receives dividends from the
company) when there are opportunities to do so. This problem might
fit better the case of R\&D investments rather than the venture
capital investments.  Namely, the company or R\&D project has a
large technology platform, based on which applications are made
and products are materialized from time to time.  Each time it
occurs, the investor tries to sell these products or applications
and in turn receives dividends. There are many papers about
dividend payout problems that consider continuous diffusion processes.
See, for example, Bayraktar and Egami \cite{BE2006c} and the
references therein.  To our knowledge, one of the few exceptions aside from \cite{kyprianou-dividend} (that we refer to earlier) is
Dassios and Embrechts \cite{dassios-embrechts} that analyze, using
the Laplace transform method, the downward jump case. The absolute
value of the jumps are exponentially distributed. In
\cite{dassios-embrechts}, the investor extracts dividends every
time when a \emph{piecewise deterministic Markov process} hits a
certain boundary (i.e., singular control). In what follows, the
dividend payments are triggered by sporadic jumps of the process
as well as the diffusion part. Whenever the value of the company
exceed a certain value, which may occur continuously or via jumps,
dividends are paid out. So the investor applies a mixture of
singular and impulse controls.

Again, we consider the jump diffusion model (\ref{eq:processjump})
for the intrinsic value of the company. Accounting for the
dividend payments the value of the company follows:
\begin{equation}\label{eq:controlled-2}
    dX_t = \mu(X_t)dt + \sigma(X_t)dW_t + \int_0^\infty y N(dt, dy)  - dZ_t
    \end{equation}
in which $Z=(Z_t)_{t\ge 0}$ is a continuous non-decreasing (expect at
$t=0$) $\mathbb{F}$-adapted process, i.e., $Z \in \mathcal{V}$, is the dividend payment policy.

The investor wants to maximize the discounted expected value of the payments she receives, which is
given by
\begin{equation} \label{eq:div_J}
J^{Z}(x)\triangleq\E^{x}\left[\int_{0}^{\tau_0}e^{-\alpha t}dZ_t
\right]
\end{equation}
in which $\tau_0$ is defined as in (\ref{eq:ruin-time})
denotes the time of insolvency. The investor wants to determine the optimal dividend policy $Z^*$
that satisfies
\begin{equation}\label{eq:div_value}
V(x) \triangleq \sup_{Z \in \mathcal{V}}J^{Z}(x)=J^{Z^*}(x),
\end{equation}
if such a $Z^* \in \mathcal{V}$ exists.
\subsection{A Mixed Singular and Impulse Control Problem}
\subsubsection{Verification Lemma}
\begin{lemma}\label{lem:harvest-verify}
Let us assume that $\sigma(\cdot)$ is bounded. If non-negative function $v \in
\mathbb{C}^1(\mathbb{R}_+)$ is also
twice continuously differentiable except at countably many points
and satisfies
\begin{enumerate}
\item [(i)] $(\A-\alpha)v(x)\leq 0$, $x \geq 0$,
\item [(ii)] $v'(x) \geq 1$, $x \geq 0$,
\item [(iii)] $v''(x) \leq 0$ (i.e. $v$ is concave),
\end{enumerate}
then
\begin{equation}\label{eq:ha-opt}
v(x) \geq V(x), \quad x \geq 0.
\end{equation}

Moreover, if there exists point $b\in \R_+$ such that $v\in
\mathrm{C}^1(\R_+)\cap \mathrm{C}^2(\R_+ \backslash\{b\})$ such
that
\begin{enumerate}
\item [(iv)] $(\A-\alpha)v(x)=0$ , $v'(x)>1$, for all $x\in [0,
b)$,
\item [(v)] $(\A-\alpha)v(x)<0$, $v(x)=x-b+v(b)$, $x>b$,
\end{enumerate}
in which the integro-differential operator $\A$ is defined by
(\ref{eq:infinite}), then
\begin{equation}
v(x)=V(x) \quad  x\in \R_+, \quad \text{and},
\end{equation}
\begin{equation}\label{eq:local-t-p-impulse}
Z_t=(X_t-b)1_{\{X_t>b\}}+L^b_t, \quad t\geq 0,
\end{equation}
in which
\begin{equation}
L_t=\int_0^{t}1_{\{X_s=b\}}dL^b_s, \quad t \geq 0,
\end{equation}
is optimal.
\end{lemma}
\begin{proof}
Let $\tau(n)$ be as in the proof of Lemma~\ref{lem:ver}. Using
It\^{o}'s formula for semimartingales (see e.g. Jacod and Shiryaev
\cite{JS})
\begin{equation}\label{eq:Ito-harvest-2}
\begin{split}
&e^{-\alpha(\tau(n)\wedge\tau_0)}v(X_{
\tau(n)\wedge\tau_0})=v(x)+\int_0^{\tau(n)\wedge\tau_0} e^{-\alpha
s}(\mathcal{A}-\alpha)v(X_s)ds-
\int_0^{\tau(n)\wedge\tau_0}e^{-\alpha s}v'(X_s)dZ^{(c)}_s
\\&+\int_0^{\tau(n)\wedge\tau_0}e^{-\alpha
s}\sigma(X_s)v'(X_s)dW_s+\int_0^{\tau(n)\wedge\tau_0}\int_0^{\infty}e^{-\alpha s}(v(X_{s-}+y)-v(X_{s-}))(N(ds,
dy)-\nu(ds,dy))
\\&+\sum_{0\leq\theta_k \leq \tau(n) \wedge
\tau_0}e^{-\alpha
\theta_k}\left(v(X_{\theta_k})-v(X_{\theta_{k-}})\right),
+\int_0^{\tau(n) \wedge \tau_0}\int_0^{\infty}e^{-\alpha s}\left(v(X_{s})-v(X_{s-}+y)\right) N(ds, dy)
\end{split}
\end{equation}
in which $\{\theta_k\}_{k \in \mathbb{N}}$ is an increasing
sequence of $\mathbb{F}$ stopping times that are the times the
process $X$ jump due to jumps in $Z$ that do not occur at the time
of Poisson arrivals. $Z^{(c)}$ is the continuous part of $Z$, i.e.,
\begin{equation}
Z_t^{(c)} \triangleq Z_t-\sum_{0 \leq s \leq t}(Z_s-Z_{s-}).
\end{equation}
The controller is allowed to choose the jump times of $Z$ to
coincide with the jump times of $N$. But this is taken into
account in (\ref{eq:Ito-harvest-2}) in the last line. Observe that
the expression on this line is zero if the jump times of $Z$ never
coincide with those of the Poisson random measure.

Equation (\ref{eq:Ito-harvest}) can be written as
\begin{equation}\label{eq:Ito-harvest}
\begin{split}
&e^{-\alpha(\tau(n)\wedge\tau_0)}v(X_{
\tau(n)\wedge\tau_0})=v(x)+\int_0^{\tau(n)\wedge\tau_0} e^{-\alpha
s}(\mathcal{A}-\alpha)v(X_s)ds-
\int_0^{\tau(n)\wedge\tau_0}e^{-\alpha s}dZ_s
\\&+\int_0^{\tau(n)\wedge\tau_0}e^{-\alpha s}(1-v'(X_s))dZ_s+\int_0^{\tau(n)\wedge\tau_0}\int_0^{\infty}e^{-\alpha s}(v(X_{s-}+y)-v(X_{s-}))(N(ds,
dy)-\nu(ds,dy))
\\&+\sum_{0\leq\theta_k \leq \tau(n) \wedge
\tau_0}e^{-\alpha
\theta_k}\left(v(X_{\theta_k})-v(X_{\theta_{k-}})+(X_{\theta_k}-X_{\theta_{k-}})v'(X_{\theta_{k-}})\right)+\int_0^{\tau(n)\wedge\tau_0}e^{-\alpha
s}\sigma(X_s)dW_s
\\&+\int_0^{\tau(n) \wedge \tau_0}\int_0^{\infty}e^{-\alpha s}\left(v(X_{s})-v(X_{s-}+y)+(y+X_{s}-X_{s-}) v'(X_{s-}+y)\right) N(ds, dy)
\end{split}
\end{equation}

After taking expectations the stochastic integral terms
vanish. Also, the concavity of $v$ implies that
\begin{equation}
v(y)-v(x)- v'(x)(y-x) \leq 0, \quad \text{for any} \quad y>x.
\end{equation}
 Now together with the, Assumptions (i), (ii), (iii) we obtain
\begin{equation*}
v(x)\geq\E^{x}\left[e^{-\alpha( \tau(n)\wedge\tau_0)}v(X_{
\tau(n)\wedge\tau_0})+\int_0^{ \tau(n)\wedge\tau_0}e^{-\alpha
s}dZ_s\right].
\end{equation*}
Equation (\ref{eq:ha-opt}) follows from the bounded and monotone
convergence theorems.

When the control $Z$ defined in (\ref{eq:local-t-p-impulse}) is
applied, the third line (\ref{eq:Ito-harvest}) is equal to
$(x-b)1_{x>b}$, since the jump times of $Z_t$ coincide with that
of the Poisson random measure $N$ except at time zero if
$X_0=x>b$. The fourth line is also zero, because $v(\cdot)$ is
linear on $[b,\infty)$. After taking expectations and then using
assumptions (iv) and (v), monotone and bounded convergence
theorems we obtain
\begin{equation}
v(x)=\ME\left[\int_0^{\tau_0}e^{-\alpha s}dZ_s\right],
\end{equation}
which proves the optimality of $Z$ and $v(\cdot)=V(\cdot)$.

\end{proof}

\subsubsection{Construction of a Candidate Solution}

As in Section~\ref{sec:cand1} we will assume that the mean measure of the Poisson random measure
$N$ is given by $\nu(dt, dy)=\lambda dt \eta e^{-\eta y}dy$, $\mu(x)=\mu$ where $\mu>0$ and
$\sigma(x)=\sigma$.

Let us define
 \begin{equation}
 v_0(x)\triangleq
B_1e^{\gamma_1 x}+B_2e^{\gamma_2 x}+B_3e^{-\gamma_3 x}, \quad x\geq 0,
\end{equation}
for $B_1$, $B_2$, $B_3 \in \mathbb{R}$ that are to be determined.
We set our candidate function to be
\begin{align}\label{eq:candidate-harv}
v(x) \triangleq
\begin{cases}
v_0(x)& x\in [0, b),\\ x-b+v_0(b), & x \in[b, \infty).
\end{cases}
\end{align}

We will choose $B_1$, $B_2$, $B_3$ and $b$ to satisfy

\begin{align}
\frac{B_1\eta}{\gamma_1-\eta}e^{\gamma_1
b}+\frac{B_2\eta}{\gamma_2-\eta}e^{-\gamma_2 b}-\frac{B_3 \eta}{\gamma_3+\eta}e^{-\gamma_3 b}+B_1e^{\gamma_1
b}+B_2e^{-\gamma_2 b}+B_3e^{-\gamma_3 b}+\frac{1}{\eta}&=0,
\label{eq:harv-1}
\\ \gamma_1B_1e^{\gamma_1 b}+\gamma_2B_2e^{\gamma_2
b}-\gamma_3B_3e^{-\gamma_3 b}&=1
\label{eq:harv-2}
\\ \gamma_1^2B_1e^{\gamma_1 b}+\gamma_2^2B_2e^{\gamma_2
b}+\gamma_3^2B_3e^{-\gamma_3 b}&=0,  \label{eq:harv-3}
\\ B_1 + B_2 +B_3&=0. \label{eq:harv-4}
\end{align}

Equation (\ref{eq:harv-1}) by explicitly evaluating
 \begin{equation*}
(\A-\alpha)v(x)=\mu v'(x)+\frac{1}{2}\sigma^2v''(x)+\lambda\left(
\int_{0}^{b-x}v(x+y)F(dy)+\int_{b-x}^{\infty}(v(b)+(x+y-b))F(y)dy\right)-(\lambda+\alpha)
v(x)
\end{equation*}
and setting it to zero. Equations (\ref{eq:harv-2}) and (\ref{eq:harv-3}) are there to enforce first and second order smooth fit at point $b$. The last equation imposes the function $v$ to be equal to zero at point zero. The vale function, $V$ satisfies this condition since whenever the value process $X$ hits level zero bankruptcy is declared.

\begin{lemma}\label{lem:har-un}
There exists unique solution $B_1, B_2, B_3$ and $b$ to the system
of equations (\ref{eq:harv-1}), (\ref{eq:harv-2}),
(\ref{eq:harv-3}), and (\ref{eq:harv-4}) if and only if the
quantity $\mu+\lambda/\eta>0$. Moreover, $B_1>0$, $B_2>0$ and
$B_3<0$.

\end{lemma}
\begin{proof}
Using  (\ref{eq:harv-1}), (\ref{eq:harv-2}), and
(\ref{eq:harv-3}), we can express $B_1$, $B_2$ and $B_3$ as
functions of $b$:  For all $b>0$, we have
\begin{align}\label{eq:harv-Ab}
B_1(b)&=\frac{e^{-\gamma_1b}}{\eta}\frac{\gamma_2\gamma_3(\eta-\gamma_1)}{\gamma_1(\gamma_2-\gamma_1)(\gamma_1+\gamma_3)}>0,\nonumber\\
B_2(b)&=\frac{e^{-\gamma_2b}}{\eta}\frac{\gamma_3\gamma_1(\gamma_2-\eta)}{\gamma_2(\gamma_2+\gamma_3)(\gamma_2-\gamma_1)}>0,\\
B_3(b)
&=-\frac{e^{\gamma_3b}}{\eta}\frac{\gamma_1\gamma_2(\eta+\gamma_3)}{\gamma_3(\gamma_1+\gamma_3)(\gamma_2+\gamma_3)}<0.\nonumber
\end{align}
Let us define
\begin{equation}
 Q(b)\triangleq B_1(b)+B_2(b)+B_3(b), \quad b \geq 0.
\end{equation}
Our claim follows once we show that the function $b \rightarrow Q(b)$, $b \geq 0$ has a unique root.
The derivative of $Q(\cdot)$
\begin{equation}
Q'(b)=B_1'(b)+B_2'(b)+B_3'(b)<0,
\end{equation}
therefore $Q(\cdot)$ is
decreasing.
Explicitly computing $Q(0)$ in (\ref{eq:harv-Ab}), we
obtain
\begin{align*}
Q(0)>0 \quad \text{if and only if} \quad
\frac{1}{\eta\gamma_1\gamma_2\gamma_3}\Big(-\gamma_1\gamma_2\gamma_3+\eta(-\gamma_1\gamma_2+\gamma_2\gamma_3+\gamma_3\gamma_1)\Big)
=\frac{\mu+\lambda/\eta}{\alpha}>0.
\end{align*}
Since
\begin{equation}
\lim_{b \rightarrow \infty}Q(b)=-\infty
\end{equation}
the claim follows.
\end{proof}
\subsubsection{Verification of Optimality}
\begin{lemma} \label{lem:harv-optimality}
Let $B_1$, $B_2$, $B_3$ and $b$ be as in Lemma~\ref{lem:har-un}.
Then $v$ defined in (\ref{eq:candidate-harv}) satisfies

(i) $(\A-\alpha)v(x) <0$ for $x\in (b, \infty)$, (ii) $v'(x)>1$ on
$x\in [0, b)$, (iii) $v''(x)< 0$ on $x\in [0, b)$.
\end{lemma}
\begin{proof}
(i): On $x\in (b, \infty)$, $v(x)=(x-b)+v_0(b)$, we compute
\begin{align*}
(\A-\alpha)v(x) &=\mu +\lambda/ \eta -\alpha (x-b)-\alpha
v_0(b^*) < \mu+\lambda/\eta -\alpha v_0(b)\\ &=\lim_{x \downarrow
b}(\A-\alpha)v(x)=\lim_{x\uparrow b}(\A-\alpha)v(x)=0.
\end{align*}
Here we used the continuity of $v(x), v'(x)$ and $v''(x)$ at
$x=b$.\\
(ii) and (iii):  Since $B_1, B_2>0$ and $B_3<0$,
\begin{equation}
v_0'''(x)=B_1 \gamma_1^3 e^{\gamma_1 x}+B_2 \gamma_2^3 e^{\gamma_2 x}
-B_3 \gamma_3^3 e^{-\gamma_3 x}>0,
\end{equation}
i.e., $v_0''(\cdot)$ is monotonically increasing in $x$.
It follows from (\ref{eq:harv-3}) that
$v''_0(b)=0$, therefore $v_0''(x)<0$ on $x\in [0,
b)$. This proves (iii).

Since $v''_0(x)<0$, $x \in \mathbb{R}_+$, $v'_0(\cdot)$
is decreasing on $ \mathbb{R}_+$.
It follows from (\ref{eq:harv-3}) that $v'_0(b)=1$. Therefore,
$v'_0(x)>1$ on $x\in [0, b)$.  This proves (ii).
\end{proof}

\begin{figure}[h]
\begin{center}
\begin{minipage}{0.3\textwidth}
\centering \includegraphics[scale=0.6]{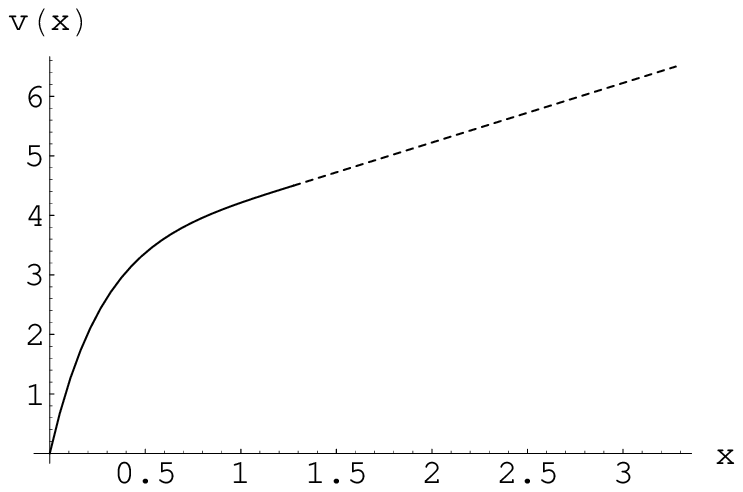} \\ (a)
\end{minipage}
\begin{minipage}{0.3\textwidth}
\centering \includegraphics[scale=0.6]{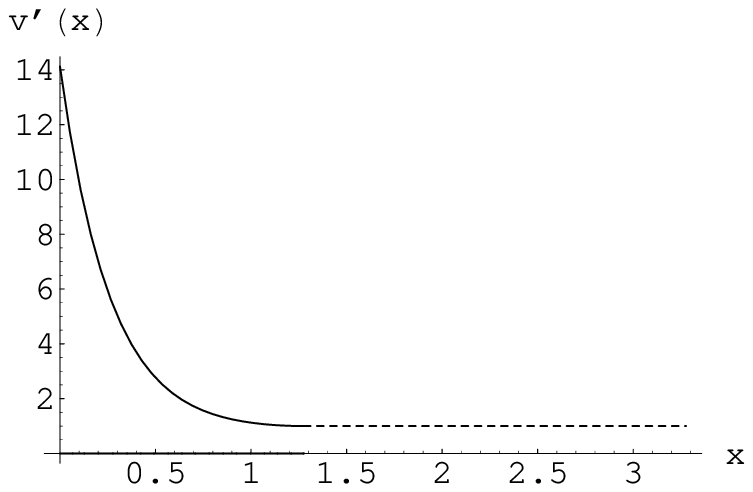} \\ (b)
\end{minipage}
\begin{minipage}{0.3\textwidth}
\centering \includegraphics[scale=0.6]{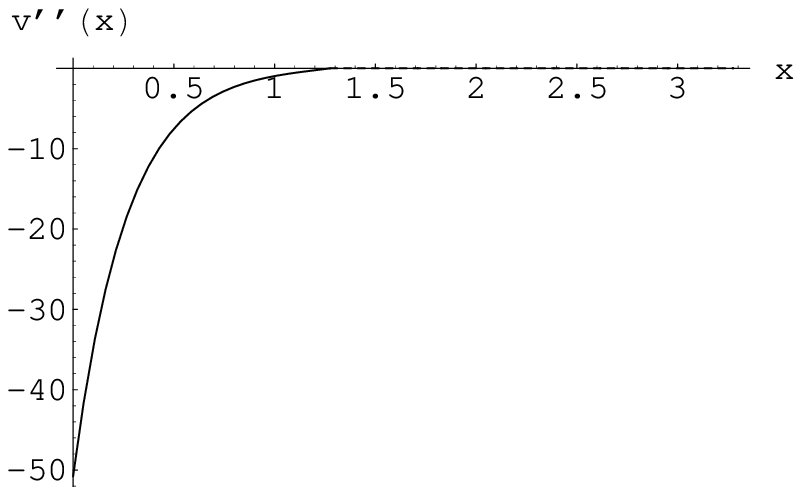} \\ (c)
\end{minipage}
\caption{\small The harvesting (dividend payout) problem with parameters
 $(\mu, \lambda, \eta, \sigma, \alpha)=(-0.05, 0.75, 1.5, 0.25, 0.1)$:  (a)
 The value function $v(x)$ with  $b=1.276$.
 (b) $v'(x)$ and $v''(x)$ to show that the optimality conditions are satisfied.}
\label{fig:harv}
\end{center}
\end{figure}

\begin{figure}[h]
\begin{center}
\begin{minipage}{0.40\textwidth}
\centering \includegraphics[scale=0.65]{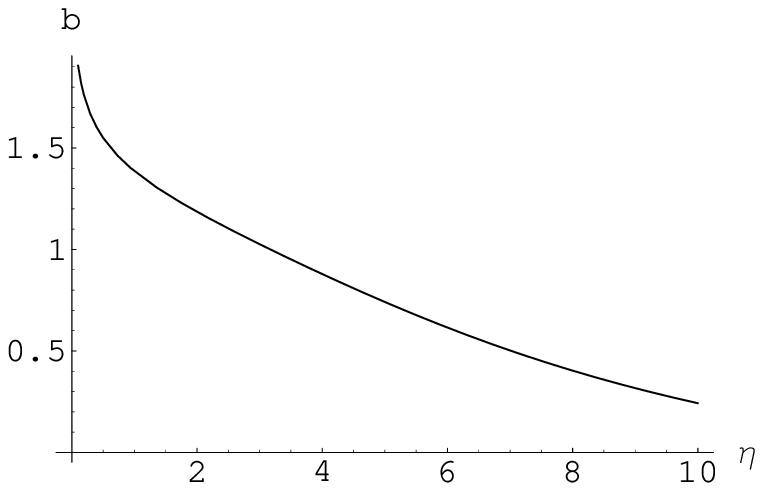} \\ (a)
\end{minipage}
\begin{minipage}{0.40\textwidth}
\centering \includegraphics[scale=0.65]{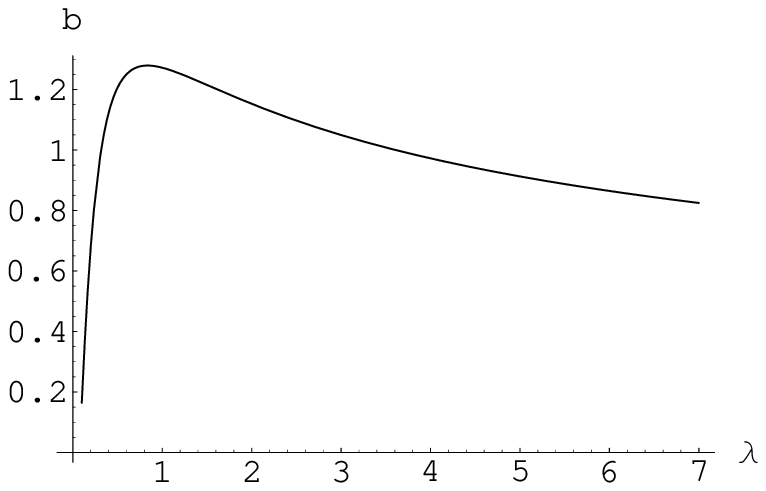} \\ (b)
\end{minipage}
\begin{minipage}{0.40\textwidth}
\centering \includegraphics[scale=0.65]{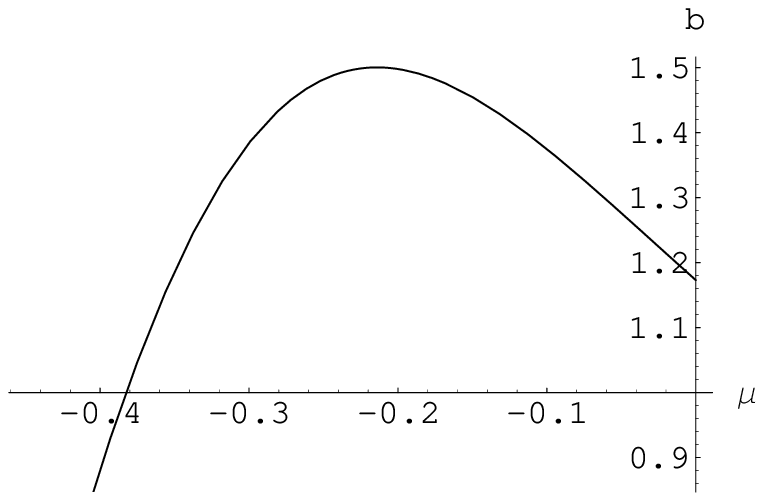} \\ (c)
\end{minipage}
\begin{minipage}{0.40\textwidth}
\centering \includegraphics[scale=0.65]{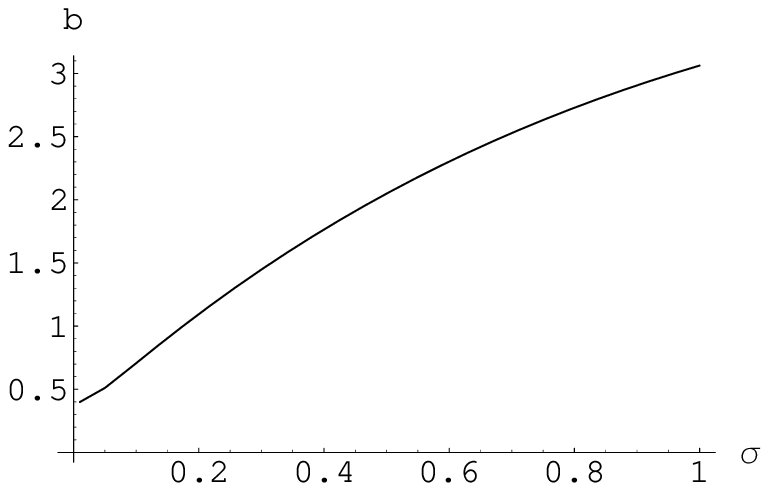} \\ (d)
\end{minipage}
\caption{\small Sensitivity analysis of the harvesting (dividend payout) problem to the parameters.  The basis parameters are $(\mu, \lambda, \eta, \sigma, \alpha)=(-0.05, 0.75, 1.5, 0.25, 0.1)$:  (a) jump  size parameter $\eta$,
(b) arrival rate $\lambda$, (c) drift rate $\mu$ and (d) volatility $\sigma$.}
 \label{fig:harv-sensitivity}
\end{center}
\end{figure}

\begin{proposition}
Suppose that $\mu+\frac{\lambda}{\eta}>0$.  Let $B_1$, $B_2$,
$B_3$ and $b$ be as in Lemma~\ref{lem:har-un}. Then the function
$v(\cdot)$ defined in (\ref{eq:candidate-harv})  satisfies
\begin{equation}
v(x)=V(x)=\sup_{Z \in \mathcal{V}}J^{Z}(x).
\end{equation}
and $Z$ defined in (\ref{eq:local-t-p-impulse}) is optimal.
\end{proposition}
\begin{proof}
Note that $(\A-\alpha)v(x)=0$, $x \in[0,b)$ as a result of (\ref{eq:harv-1}). The function $v(\cdot)$ is linear on $[b,\infty)$. It follows from Lemma~\ref{lem:harv-optimality} that the function $v(\cdot)$ satisfies all the conditions in the verification lemma.
\end{proof}
Figure \ref{fig:harv} shows the value function and its derivatives. As expected the value function is concave and is twice continuously differentiable.
Finally, we perform some sensitivity analysis of the optimal barrier $b$ with respect to the parameters of the problem.  Figure \ref{fig:harv-sensitivity} shows the results.  \textbf{Graph (a):  }  The first graph shows that as the expected value of jump size $1/\eta$ decreases, so does the threshold level $b$, as one would expect.   \textbf{Graph (b):  }  It is interesting to observe that $b^*$ increases first and start decreasing when $\lambda$ reaches a certain level, say $\lambda_{\max}$.  A possible interpretation is as follows:  In the range of $(0, \lambda_{\max})$, i.e. for small $\lambda$, one wants to extract a large amount of cash whenever jumps occur since the opportunities are limited.  As $\lambda$ gets larger, one starts to be willing to let the process live longer by extracting smaller amounts each time.
On the other hand, after $\lambda \geq \lambda_{\max}$, one becomes comfortable with receiving more dividends, causing the declining trend of $b^*$.    \textbf{Graph (c):  }  The small $\mu$ in the absolute value sense implies that it takes more time to hit the absorbing state.  Accordingly, it is safe to extract a large amount of dividend.  However, when the cost increases up to a certain level, say $\mu^*$, it becomes risky to extract and hence $b^*$ increases.  It is observed that after the cost level is beyond $\mu^*$, one would become more desperate to take a large dividend at one time in the fear of imminent insolvency caused by a large $\mu$ (in the absolute value sense).  This is the downward trend of $b^*$ on the left side of $\mu^*$.  \textbf{Graph (d):  }  As the volatility goes up, then the process tends to spend more time away from zero in both the positive and negative real line.  Accordingly, the threshold level increases to follow
 the process.
\section{Concluding Remarks}
\begin{figure}[h]
\begin{center}
\begin{minipage}{0.3\textwidth}
\centering \includegraphics[scale=0.6]{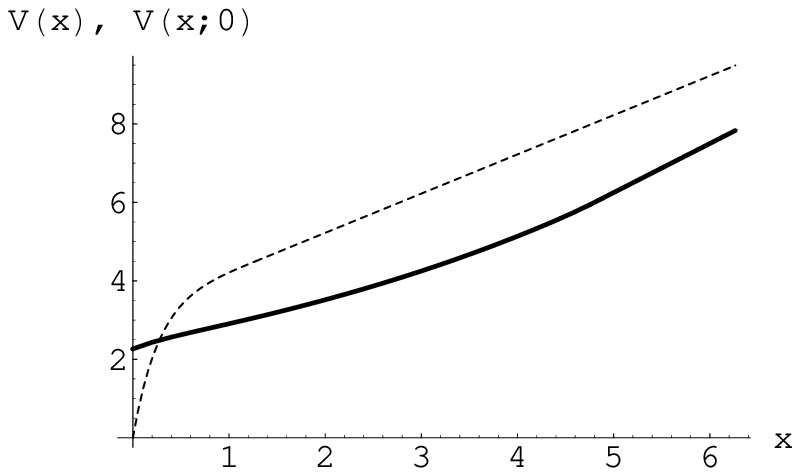} \\ (a)
\end{minipage}
\begin{minipage}{0.3\textwidth}
\centering \includegraphics[scale=0.6]{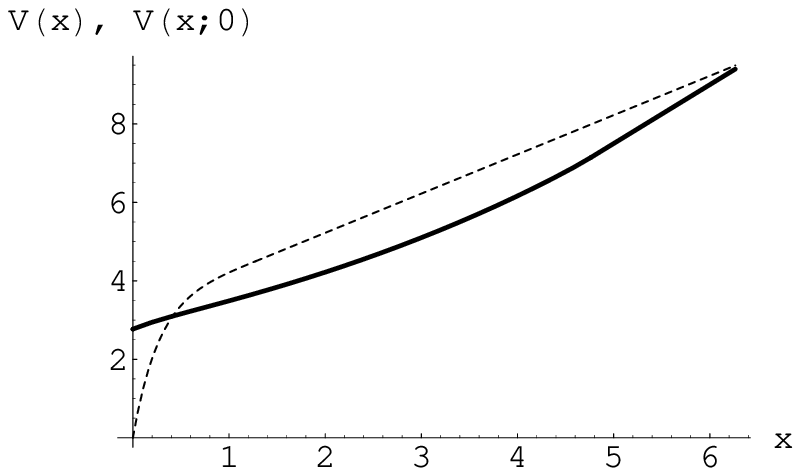} \\ (b)
\end{minipage}
\begin{minipage}{0.3\textwidth}
\centering \includegraphics[scale=0.6]{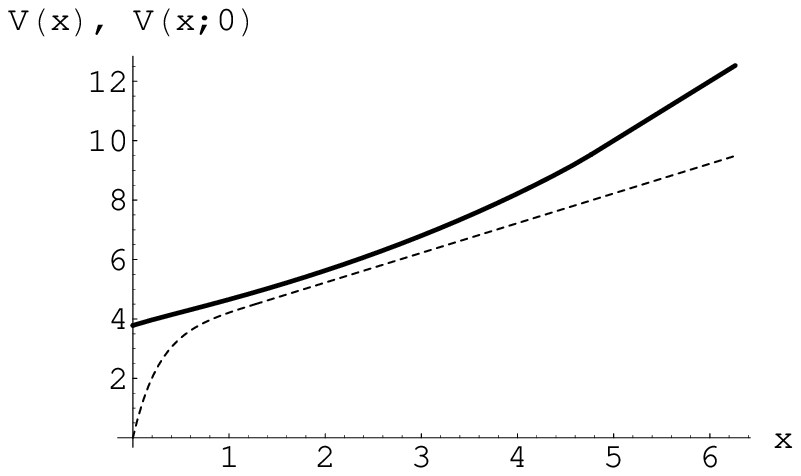} \\ (c)
\end{minipage}
\caption{\small The comparison of two value functions with
 $(\mu, \lambda, \eta, \sigma, \alpha)=(-0.05, 0.75, 1.5, 0.25, 0.1)$:  (a) $r=1.25$.
 (b) $r=1.5$ and (c) $r=2$ where the value function for the IPO problem with $a=0$ is shown in solid line and
 the value function for the harvesting problem is shown in dashed line.}
\label{fig:comparison}
\end{center}
\end{figure}

Before concluding, we compare two value functions, one for the IPO
problem and the other for the harvesting problem.  We set
parameters $\mu, \sigma, \lambda, \eta$, and $\alpha$ equal and
vary the level of $r>1$, the expected return at the IPO market.
Figure \ref{fig:comparison} exhibits the two value functions:
$v(x; 0)$ (solid line) for the IPO problem with $a=0$ and $v(x)$
(dashed line) for the harvesting problem.  We consider three
different values of $r$ here; (a) $r=1.25$, (b) $r=1.5$ and (c)
$r=2$.  It can be observed that as $r$ increases, the value
function for the IPO problem shifts upward for all the points of
$x \in\R_+$.  This jump diffusion model, although simple, gives a
quick indication as to which strategy (IPO or harvesting) is more
advantageous given the initial investment amount $x$.  Moreover,
as we discussed, this model has both diffusion and jump
components, allowing us to model different stage of the start-up
company by modifying the relative size of diffusion parameter
$\sigma$ and jump parameter $\lambda/\eta$.

\end{document}